\documentclass[11pt, oneside, reqno]{amsart}
\usepackage[a4paper,margin=1in]{geometry}
\usepackage[T1]{fontenc}
\usepackage{amsfonts}
\usepackage{amssymb}
\usepackage{amsmath}
\usepackage{amsthm}
\usepackage{mathtools}
\usepackage{mathrsfs}
\usepackage{zref-clever}
\usepackage{hyperref}
\usepackage{stmaryrd}
\usepackage{tikz-cd}
\usepackage{enumitem}
\usepackage{multirow}
\usepackage{bbm}
\usepackage{comment}
\usepackage{float}
\usepackage{url}
\usepackage[en-GB]{datetime2}
\DTMlangsetup[en-GB]{ord=omit}
\allowdisplaybreaks

\numberwithin{equation}{section}

\zcsetup{cap,nameinlink=false}

\AddToHook{env/lemma/begin}{%
	\zcsetup{countertype={theorem=lemma}}}
\AddToHook{env/proposition/begin}{%
	\zcsetup{countertype={theorem=proposition}}}
\AddToHook{env/corollary/begin}{%
	\zcsetup{countertype={theorem=corollary}}}
\AddToHook{env/conjecture/begin}{%
	\zcsetup{countertype={theorem=conjecture}}}
\AddToHook{env/remark/begin}{%
	\zcsetup{countertype={theorem=remark}}}
\AddToHook{env/example/begin}{%
	\zcsetup{countertype={theorem=example}}}

\zcRefTypeSetup{question}{
	Name-sg = Question ,
	name-sg = question ,
	Name-pl = Questions ,
	name-pl = questions ,
}

\zcRefTypeSetup{conjecture}{
	Name-sg = Conjecture ,
	name-sg = conjecture ,
	Name-pl = Conjectures ,
	name-pl = conjectures ,
}

\zcRefTypeSetup{observation}{
	Name-sg = Observation ,
	name-sg = observation ,
	Name-pl = Observations ,
	name-pl = observations ,
}

\zcRefTypeSetup{speculation}{
	Name-sg = Speculation ,
	name-sg = speculation ,
	Name-pl = Speculations ,
	name-pl = speculations ,
}

\zcRefTypeSetup{remark}{
	Name-sg = Remark ,
	name-sg = remark ,
	Name-pl = Remarks ,
	name-pl = remarks ,
}

\zcRefTypeSetup{example}{
	Name-sg = Example ,
	name-sg = example ,
	Name-pl = Examples ,
	name-pl = examples ,
}

\newtheorem{theorem}{Theorem}[section]
\newtheorem{lemma}[theorem]{Lemma}
\newtheorem{proposition}[theorem]{Proposition}
\newtheorem{corollary}[theorem]{Corollary}

\theoremstyle{definition}

\theoremstyle{remark}
\newtheorem{remark}[theorem]{Remark}

\def\MR#1{}
\newcommand\m[1]{\begin{pmatrix}#1\end{pmatrix}}
\newcommand\bm[1]{\begin{bmatrix}#1\end{bmatrix}}
\newcommand\sm[1]{\begin{psmallmatrix}#1\end{psmallmatrix}}
\newcommand\bsm[1]{\begin{bsmallmatrix}#1\end{bsmallmatrix}}

\newcommand\cl[1]{\mathcal{#1}}

\newcommand\wt[1]{\widetilde{#1}}
\newcommand\wh[1]{\widehat{#1}}

\newcommand{\Z}{\mathbb{Z}}
\newcommand{\R}{\mathbb{R}}
\newcommand{\C}{\mathbb{C}}

\newcommand{\F}{\mathbb{F}}

\newcommand{\clh}{\mathcal{H}}

\newcommand{\eq}{\;=\;}
\newcommand{\deq}{\;:=\;}

\def\thin{{\hspace{1pt}}}

\newcommand{\id}{\mathrm{id}}

\newcommand{\M}{\mathsf{M}}

\newcommand{\GL}{\mathrm{GL}}

\newcommand{\PSL}{\mathrm{PSL}}
\newcommand{\SL}{\mathrm{SL}}

\newcommand{\Ker}{\mathrm{Ker}}

\newcommand{\sgn}{\mathrm{sgn}}

\newcommand{\Span}{\mathrm{Span}}

\renewcommand{\Re}{\operatorname{Re}}

\newmuskip\pFqskip
\mathchardef\pFcomma=\mathcode`,

\title{Twisting-free converse theorems for $\Gamma_0(32)$ and $\Gamma_0(36)$}
\author{Pengcheng Zhang}
\address{}
\email{}
\date{\today}

\begin{document}

\begin{abstract}
    We prove twisting-free converse theorems for $\Gamma_0(32)$ and $\Gamma_0(36)$ under the setting of Conrey--Farmer, namely assuming a single functional equation and a suitable Euler product, with particular shapes at ramified primes. The proofs make substantial use of the ramified primes, with methods somewhat different from the existing ones in the literature.
\end{abstract}

\maketitle

\section{Introduction}

Converse theorems for modular forms date back to Hecke \cite{hecke}, who settled the result for level $1$ assuming a single functional equation. Weil \cite{weil} generalized this to a converse theorem for level $N$ by assuming all the functional equations twisted by Dirichlet characters. Later in \cite{conrey-farmer}, Conrey and Farmer proved certain \emph{twisting-free} converse theorems for a few small levels by assuming only a single functional equation together with an Euler product. Their result was further generalized to some other levels by \cite{conrey-farmer-odgers-snaith,bedert-cooper-oliver-zhang}.

Following the setting of Conrey--Farmer, we will prove twisting-free converse theorems for $\Gamma_0(32)$ and $\Gamma_0(36)$ in this paper. We here recall the setting and state the main results.

\subsection{Main results}

Let $N\in\Z^+$, $k\in2\Z^+$, and $(a_n)$ be a sequence of complex numbers of polynomial growth. Let
\begin{align*}
	f(z)\eq\sum_{n=1}^\infty a_ne^{2\pi inz}
\end{align*}
be a function on the upper half plane. Define an $L$-series
\begin{align*}
	L(s,f)\eq\sum_{n=1}^\infty\frac{a_n}{n^s}
\end{align*}
and its completed version
\begin{align*}
	\widehat{L}(s,f)\eq N^{s/2}(2\pi)^{-s}\Gamma(s)L(s,f)
\end{align*}
for $\Re(s)\gg 0$.

\begin{theorem}
\zlabel{mainthm}
Let $N=32$ or $36$, $k\in 2\Z^+$, and $\epsilon=\pm1$. Suppose that
	\begin{enumerate}[label=\rm{(\alph*)}]
		\item $\wh{L}(s,f)$ extends to an entire function that is bounded in vertical strips and satisfies a functional equation
		\begin{align*}
			\wh{L}(s,f)\eq\epsilon\thin i^k\thin\wh{L}(k-s,f);
		\end{align*}
		\item $L(s,f)$ has an Euler product of the form
		\begin{align*}
			L(s,f)\eq\prod_p L_p(s,f),
		\end{align*}
		where
        \begin{align}
        \label{local-euler}
            L_p(s,f)\eq\begin{cases}(1-a_pp^{-s}+p^{k-1-2s})^{-1} &\text{ if } p\nmid N, \\
            1 &\text{ if } p^2\mid N.
            \end{cases}
        \end{align}
	\end{enumerate}
	Then, $f$ is a cusp form of level $N$ and weight $k$.
\end{theorem}

\begin{remark}
    Conrey and Farmer also imposed the condition
    \begin{align}
    \label{local-euler-pmmn}
        L_p(s,f)\eq (1-p^{k/2-1-s})^{-1} \text{ if } p \mathrel\Vert N
    \end{align}
    for the Euler product, which is irrelevant for the levels considered here. Assuming both (\ref{local-euler}) and (\ref{local-euler-pmmn}), the modularity of $f$ was known for $1\leq N\leq 18$ and $N=20,23,24$ by \cite{hecke}, \cite{conrey-farmer}, \cite{conrey-farmer-odgers-snaith}, and \cite[Theorem~2.4]{bedert-cooper-oliver-zhang}. Among these cases, cuspidality was known for $1\leq N\leq 17$ and $N=23$ by \cite{conrey-farmer,conrey-farmer-odgers-snaith}. We also refer to \cite[Theorem~2.1]{conjhecke} and \cite[Theorem~7.2]{zhang-thesis} for results involving assumptions only on unramified primes.
\end{remark}

\begin{remark}
    One can relax the full Euler product in (b) to a partial one. That is, one can choose a finite set $S$ of primes and consider a partial Euler product of the form
        \begin{align*}
            L(s,f)\eq\prod_{p\in S}L_p(s,f)\cdot\bigg(\sum_{(n,S)=1}\frac{a_n}{n^s}\bigg),
        \end{align*}
        where $L_p(s,f)$ is defined the same as before and the sum takes over all positive integers $n$ that are coprime to primes in $S$. Then, the theorem holds for $(N,S)=(32,\{2,3\})$ and $(36,\{2,3,5\})$.
\end{remark}

We briefly describe the proof strategy here. In both cases, the idea is to consider $h=f|_k\sm{1&0\\0&d}$, where $d=2$ for $N=32$ and $d=6$ for $N=36$. We will prove that $h$ is a cusp form of level $\Gamma_0(8)$ and $\Gamma(6)$ respectively, from which the modularity and cuspidality of $f$ follow by conjugation.

\subsection{Acknowledgements}

This research was partly conducted during the author's stay at Max Planck Institute for Mathematics. The author therefore thanks Max Planck Institute for Mathematics for their hospitality and relevant support.

\subsection{Disclosure of AI use}

The proofs of the two converse theorems were originally generated by GPT. These proofs are available at \url{https://github.com/chbe113/ai-related/tree/main/ConverseThm32and36}, which may be of independent interest. The author digested both proofs, partly with the assistance of GPT, and incorporated them into this paper according to his own understanding. This paper itself was written solely by the author, and the author takes full responsibility for the correctness of this paper.

\section{Preliminary lemmas}

We first set up some notation and recall several lemmas from \cite{conrey-farmer}. For $k\in\Z$, $\gamma=\sm{a&b\\c&d}\in\GL_2^+(\R)$, and $g:\clh\rightarrow\C$, define the \emph{slash operator of weight $k$} by
\begin{align*}
    (g|_k\gamma)(z)\deq (\det\gamma)^{k/2}(cz+d)^{-k}g\bigg(\frac{az+b}{cz+d}\bigg).
\end{align*}
In our case, since $k\in 2\Z^+$, any scalar matrix $\lambda I$ for nonzero $\lambda\in\R$ always acts trivially. Thus, we will typically regard $\gamma\in\GL_2^+(\R)$ as an element in $\SL_2(\R)$ up to rescaling. We will also extend the slash operator linearly to the group ring $\C[\SL_2(\R)]$. To avoid possible confusion, we write $\gamma$ as an element in the matrix group and $[\gamma]$ in the group ring.

For $a\in\R$ and $M\in\Z$, let
\begin{align*}
    T^a&\eq\m{1&a\\0&1}, \\
    W_M&\eq\m{0&-1\\M&0}, \\
    A_M&\eq W_MT^{-1}W_M^{-1}\eq \m{1&0\\M&1}.
\end{align*}
For a prime $p$, let
\begin{align*}
    \cl{T}_p&\eq\bm{p&0\\0&1}+\sum_{a=0}^{p-1}\bm{1&a\\0&p}, \\
    \cl{U}_p&\eq\sum_{a=0}^{p-1}\bm{1&a\\0&p}.
\end{align*}

\begin{lemma}
\zlabel{preliminary}
Under the assumptions in \zcref{mainthm}, we have
\begin{enumerate}[label=\rm{(\arabic*)}]
    \item $f|_kT=f$;
    \item $f|_kW_N=\epsilon f$;
    \item for any prime $p$,
    \begin{alignat*}{2}
    f|_k\cl{T}_p&=p^{1-k/2} a_pf && \quad\text{ if }  p \nmid N,\\
    f|_k\thin\cl{U}_p&=0 && \quad\text{ if }  p^2 \mid N.
    \end{alignat*}
\end{enumerate}
\end{lemma}
\begin{proof}
    This follows from standard computations. See \cite[Lemma~1]{conrey-farmer} (and also \cite[Proposition~11.3.1]{cohen-stromberg}).
\end{proof}

\begin{corollary}
\zlabel{4midN}
    Under the assumptions in \zcref{mainthm}, if $4\mid N$, then we have
    \begin{align*}
        f|_kT^{1/2}\eq f|_k A_{N/2}\eq -f.
    \end{align*}
\end{corollary}
\begin{proof}
   As $4\mid N$, we have $f|_k\big(\bsm{1&0\\0&2}+\bsm{1&1\\0&2}\big)=0$, so $f|_kT^{1/2}=f|_k\sm{1&1\\0&2}\sm{1&0\\0&2}^{-1}=-f$. The anti-invariance under $A_{N/2}$ now follows from \zcref{preliminary} and the fact that $A_{N/2}=W_NT^{-1/2}W_N^{-1}$.
\end{proof}

\begin{proposition}
\zlabel{general}
Let $p$ be a prime with $p\nmid N$. Under the assumptions in \zcref{mainthm}, we have
\begin{align*}
    f|_k\sum_a\bigg(\bm{p&-a \\ -bN&(1+abN)/p}-1\bigg)[T^{a/p}]\eq 0,
\end{align*}
where the sum is taken over a set of representatives of $(\Z/p\Z)^\times$, and $b$ is chosen such that $p\mid 1+abN$.
\end{proposition}
\begin{proof}
    See \cite[Corollary~2]{conrey-farmer}.
\end{proof}

The next lemma is a result of the holomorphic condition. This lemma plays a key role in Ogg's proof \cite[Chapter~5]{ogg} of Weil's converse theorem as well as the proof of the converse theorem by Conrey--Farmer. 

\begin{lemma}
	\zlabel{weil}
	Let $k\in\Z$ with $k\neq 0$. Suppose that $g:\clh\rightarrow\C$ is holomorphic and $e\in\SL_2(\R)$ is an elliptic matrix of infinite order. Then, $g|_ke=g$ implies that $g=0$. 
\end{lemma}
\begin{proof}
    See \cite[Lemma~11.3.5]{cohen-stromberg}.
\end{proof}

We finish this section with an easy lemma that will be needed for establishing the cuspidality.

\begin{lemma}
\zlabel{cusp-lemma}
    Let $g\in M_k(\Gamma(N))$. Let $\gamma\in\SL_2(\Z)$ and suppose that there exists $w>0$ such that $(g|_k\gamma)(z+w)=-(g|_k\gamma)(z)$. Then, $g|_k\gamma$ vanishes at infinity.
\end{lemma}
\begin{proof}
    Since $(g|_k\gamma)(z+2w)=(g|_k\gamma)(z)$, we have a Fourier series expansion
    \begin{align*}
        (g|_k\gamma)(z)\eq\sum_{n\in\Z}^\infty b_ne^{\pi inz/w}.
    \end{align*}
    As $g$ is holomorphic at every cusp, $g|_k\gamma$ is also holomorphic at every cusp, so $b_n=0$ for $n<0$. Now, the assumption that $(g|_k\gamma)(z+w)=-(g|_k\gamma)(z)$ implies that $b_0=0$. The result then follows.
\end{proof}

\section{Proof of \zcref{mainthm}}

We will now prove \zcref{mainthm}. Throughout this section, we adopt all the assumptions in \zcref{mainthm}.

\subsection{The case when $N=32$}

We first prove \zcref{mainthm} for $N=32$. As $4\mid 32$, we have
\begin{align*}
    f|_kA_{16}\eq f|_kT^{1/2}\eq -f
\end{align*}
by \zcref{4midN}.

\begin{lemma}
\zlabel{32-f-inv-mat}
We have
\begin{align*}
    f|_k\m{3&1\\32&11}\eq f.
\end{align*}
\end{lemma}
\begin{proof}
Applying \zcref{general} with $p=3$ yields that
    \begin{align*}
        f|_k\bigg(\bm{3&-1\\-32&11}-1\bigg)[T^{1/3}]+f|_k\bigg(\bm{3&1\\32&11}-1\bigg)[T^{-1/3}]\eq 0.
    \end{align*}
    Note that
    \begin{align*}
        \m{3&-11\\-32&11}\eq A_{16}^{-1}T^{1/2}\m{3&1\\32&11}^{-1}A_{16}T^{-1/2}
    \end{align*}
    and that $f|_kA_{16}=f|_kT^{1/2}=-f$. Write $g=f|_k\big(\bsm{3&1\\32&11}-1\big)$. We then obtain that
    \begin{align*}
        g\eq  g|_k \m{3&1\\32&11}^{-1}A_{16}T^{-1/2}T^{2/3},
    \end{align*}
    where
    \begin{align*}
        \m{3&1\\32&11}^{-1}A_{16}T^{-1/2}T^{2/3}\eq\m{-5&-11/6\\16&17/3} 
    \end{align*}
    is elliptic of infinite order. Thus, $g=0$ by \zcref{weil}, and the result follows.
\end{proof}

Let
\begin{align*}
    h\eq f|_k\m{1&0\\0&2}.
\end{align*}

\begin{lemma}
\zlabel{32-h-transform}
We have
\begin{align*}
    h|_kT\eq h|_kA_8\eq -h,\qquad h|_k\m{3&2\\16&11}\eq h.
\end{align*}
\end{lemma}
\begin{proof}
    These simply follow from the identities
    \begin{align*}
    f|_kT^{1/2}&\eq f|_kA_{16}\eq -f,\qquad f|_k\m{3&1\\32&11}\eq f. \qedhere
    \end{align*}
\end{proof}

\begin{lemma}
\zlabel{gamma08-generator}
We have
\begin{align*}
    \Gamma_0(8)\eq\bigg\langle T,\,A_8,\,\m{3&2\\16&11},\,-I\bigg\rangle.
\end{align*}
\end{lemma}
\begin{proof}
    By SageMath \cite{sagemath}, we have
    \begin{align*}
        \Gamma_0(8)\eq\bigg\langle\m{1&1\\0&1},\m{5&-1\\16&-3},\m{5&-2\\8&-3},\m{-1&0\\0&-1}\bigg\rangle.
    \end{align*}
    We also have the following equalities.
    \begin{align*}
        \m{5&-1\\16&-3}\eq -T\m{3&2\\16&11}^{-1},\quad \m{5&-2\\8&-3}\eq -TA_8^{-1}\m{3&2\\16&11}T^{-1}.
    \end{align*}
    The result then follows.
\end{proof}

Now, define $\chi:\Gamma_0(8)\rightarrow\{\pm1\}$ by setting
\begin{align*}
    \chi(\gamma)\eq(-1)^{b+c}
\end{align*}
for $\gamma=\sm{a&b\\8c&d}\in\Gamma_0(8)$. It is easy to check that $\chi$ defines a character and that $\chi(T)=\chi(A_8)=-1$ and $\chi\big(\sm{3&2\\16&11}\big)=1$.

\begin{corollary}
\zlabel{h-gamma08-inv}
    We have $h\in S_k(\Gamma_0(8),\chi)$.
\end{corollary}
\begin{proof}
    It follows from \zcref{32-h-transform,gamma08-generator} that $h|_k\gamma=\chi(\gamma)h$ for all $\gamma\in\Gamma_0(8)$. Now, as $h=f|_k\sm{1&0\\0&2}$, the definition of $f$ implies that $h$ has a Fourier expansion with only positive terms and coefficients of polynomial growth. By \cite[Proposition~1.2.4]{diamond-shurman}, this implies that $h\in M_k(\Gamma_0(8),\chi)$.

    Now, we prove the cuspidality. For a cusp $c$ of $\Gamma_0(8)$ with width $w$, one chooses some $\alpha\in\SL_2(\Z)$ with $\alpha\infty=c$ and computes the stabilizer $P=\alpha T^w\alpha^{-1}$. Performing this for all the four cusps of $\Gamma_0(8)$ gives the following table.
    \[
\begin{array}{c|c|c}
\text{Cusp $c$} & \text{Width $w$} & \text{Stabilizer $P$}  \\ \hline
\infty & 1 &
\sm{1&1\\0&1} \\[6pt]
0 & 8 &
\sm{1&0\\-8&1} \\[6pt]
\frac12 & 2 &
\sm{-3&2\\-8&5} \\[6pt]
\frac14 & 1 &
\sm{-3&1\\-16&5}
\end{array}
\]
It is easy to see that $\chi(P)=-1$ for any cusp $c$. Hence, $h|_kP=-h$. In particular, this implies that $(h|_k\alpha)(z+w)=-(h|_k\alpha)(z)$. As $h\in M_k(\Gamma_0(8),\chi)$, by \zcref{cusp-lemma}, this implies that $h|_k\alpha$ vanishes at infinity, so $h$ vanishes at the corresponding cusp. In conclusion, $h$ vanishes at all the cusps of $\Gamma_0(8)$, so it is a cusp form.
\end{proof}

\begin{corollary}
    \zcref{mainthm} holds when $N=32$. That is, $f\in S_k(\Gamma_0(32))$.
\end{corollary}
\begin{proof}
    For $\sm{a&b\\32c&d}\in\Gamma_0(32)$, we have
    \begin{align*}
        \m{1&0\\0&1/2}\m{a&b\\32c&d}\m{1&0\\0&2}\eq\m{a&2b\\16c&d}\in\Gamma_0(8)
    \end{align*}
    with $\chi\big(\sm{a&2b\\16c&d}\big)=(-1)^{2b+2c}=1$, so $h|_k\sm{a&2b\\16c&d}=h$ by \zcref{h-gamma08-inv}. As $f=h|_k\sm{1&0\\0&1/2}$ and $h\in S_k(\Gamma_0(8),\chi)$ by \zcref{h-gamma08-inv}, it follows that $f\in S_k(\Gamma_0(32))$.
\end{proof}

    In fact, one could already deduce the invariance of $f$ under $\Gamma_0(32)$ after obtaining that $f|_k\sm{3&1\\32&11}=f$ by \zcref{32-f-inv-mat}. Indeed, by SageMath \cite{sagemath}, we have
    \begin{align*}
        \Gamma_0(32)\eq\bigg\langle &\begin{pmatrix}1&1\\0&1\end{pmatrix},
\begin{pmatrix}9&-1\\64&-7\end{pmatrix},
\begin{pmatrix}19&-3\\32&-5\end{pmatrix},
\begin{pmatrix}79&-14\\96&-17\end{pmatrix},
\begin{pmatrix}49&-9\\256&-47\end{pmatrix}, \\
&\begin{pmatrix}9&-2\\32&-7\end{pmatrix},
\begin{pmatrix}25&-9\\64&-23\end{pmatrix},
\begin{pmatrix}27&-11\\32&-13\end{pmatrix},
\begin{pmatrix}25&-18\\32&-23\end{pmatrix},
\begin{pmatrix}-1&0\\0&-1\end{pmatrix}\bigg\rangle.
\end{align*}
   Writing $\gamma=\sm{3&1\\32&11}$, we have the following equalities.
\begin{alignat*}{2} 
\begin{pmatrix}9&-1\\64&-7\end{pmatrix}
&\eq (T^{1/2}\gamma^{-1})^2,&\quad
\begin{pmatrix}19&-3\\32&-5\end{pmatrix}
&\eq T^{1/2}\gamma T^{-1/2},\\
\begin{pmatrix}79&-14\\96&-17\end{pmatrix}
&\eq -T(A_{16}^{-1}\gamma T^{-1/2})^2,&\quad 
\begin{pmatrix}49&-9\\256&-47\end{pmatrix}
&\eq T^{1/2}\gamma^{-1}A_{16}T^{-1}
  A_{16}^{-1}\gamma T^{-1/2},\\
\begin{pmatrix}9&-2\\32&-7\end{pmatrix}
&\eq T^{1/2}(A_{16}^{-1}\gamma)^2T^{-1/2},&\quad
\begin{pmatrix}25&-9\\64&-23\end{pmatrix}
&\eq T^{1/2}A_{16}^{-1}
  (T^{1/2}\gamma^{-1})^2A_{16}T^{-1/2},\\
\begin{pmatrix}27&-11\\32&-13\end{pmatrix}
&\eq TA_{16}^{-1}\gamma T^{-1/2}A_{16}T^{-1/2},&\quad
\begin{pmatrix}25&-18\\32&-23\end{pmatrix}
&\eq T(A_{16}^{-1}\gamma)^2T^{-1}.
\end{alignat*}
It then follows that $f|_k\gamma=f$ for all $\gamma\in\Gamma_0(32)$.

\subsection{The case when $N=36$}

We now prove \zcref{mainthm} for $N=36$. As $4\mid 36$ and $9\mid 36$, we have
\begin{align*}
    f|_kT^{1/2}\eq -f,\qquad f|_k(1+[T^{1/3}]+[T^{2/3}])\eq0
\end{align*}
by \zcref{preliminary}. The two equalities together also imply that
\begin{align*}
    f|_k(1-[T^{1/6}]+[T^{1/3}])\eq 0.
\end{align*}
Now, let
\begin{align*}
    h\eq f|_k\m{1&0\\0&6}.
\end{align*}
In the following part, we also adopt the standard notation
\begin{align*}
    S\eq\m{0&-1\\1&0}.
\end{align*}

\begin{lemma}
\zlabel{36-h-transform}
We have
\begin{align*}
    h|_k(1-[T]+[T^2])\eq 0,\qquad h|_kT^3\eq-h,\qquad h|_kS\eq\epsilon h.
\end{align*}
\end{lemma}
\begin{proof}
    These simply follow from the following equalities
    \begin{align*}
        f|_k(1-[T^{1/6}]+[T^{1/3}])&\eq 0,\qquad f|_kT^{1/2}\eq-f,\qquad f|_kW_{36}\eq\epsilon f. \qedhere
    \end{align*}
\end{proof}

Now, let
\begin{align*}
    V\eq\Span_\C\{ h|_k\gamma\mid \gamma\in\SL_2(\Z)\}
\end{align*}
and let
\begin{align*}
    h_1\eq h|_kT,\qquad h_2\eq h|_kTS\eq h_1|_kS.
\end{align*}
By \zcref{36-h-transform}, we obtain the following table.
\begin{equation}
    \label{table}
    \begin{array}{c|cc}
    & |_kS & |_kT\\ \hline
    h & \epsilon h & h_1\\
    h_1 & h_2 & h_1-h\\
    h_2 & h_1 & h-\epsilon h_2
    \end{array}
\end{equation}
Here only the bottom right entry is not a straightforward computation. For this one, we have
\begin{align*}
    h_2|_kT\eq h|_kTST\eq h|_kST^{-1}S\eq \epsilon h|_k T^{-1}S\eq \epsilon (h|_kS-h|_kTS)\eq h-\epsilon h_2,
\end{align*}
where we use the fact that $TST=-ST^{-1}S$ and that $h|_k([T^{-1}]-1+[T])=0$. Equivalently, we can write (\ref{table}) as the following
\begin{align}
\label{st-matrices}
\begin{aligned}
    \m{h&h_1&h_2}|_kT&\eq\m{h&h_1&h_2}\m{0&-1&1\\1&1&0\\0&0&-\epsilon}, \\
    \m{h&h_1&h_2}|_kS&\eq\m{h&h_1&h_2}\m{\epsilon &0&0\\0&0&1\\0&1&0}.
\end{aligned}
\end{align}
It is easy to check that $\sm{0&-1&1\\1&1&0\\0&0&-\epsilon}^6=I$ for $\epsilon=\pm1$. Hence, $T^6$ acts trivially on $h$, $h_1$, and $h_2$.

\begin{corollary}
\zlabel{shape-of-v}
We have
\begin{align*}
    V\eq\Span_\C\{h,h_1,h_2\}.
\end{align*}
In particular, $V$ is finite-dimensional and $T^6$ acts trivially on $V$.
\end{corollary}
\begin{proof}
    Since $\SL_2(\Z)$ is generated by $S$ and $T$, by (\ref{table}), the space $\Span_\C\{h,h_1,h_2\}$ is closed under the action of $\SL_2(\Z)$, and $h|_k\gamma\in\Span_\C\{h,h_1,h_2\}$ for all $\gamma\in\SL_2(\Z)$. The result then follows.
\end{proof}

To proceed with the proof, we will show that $T^3$ acts as $-\id_V$ on $V$. For this, we need to use the Euler product at $5$.

\begin{lemma}
\zlabel{t3-triv}
    The matrix $T^3$ acts as $-\id_V$ on $V$.
\end{lemma}
\begin{proof}
    If $\epsilon=1$, then it is easy to check that $\sm{0&-1&1\\1&1&0\\0&0&-\epsilon}^3=-I$, so the result follows from (\ref{st-matrices}). Hence, suppose that $\epsilon=-1$. 
    
    Applying \zcref{general} with $p=5$ yields that
\begin{align*}
    f|_k\bigg(&\bigg(\bm{5&-1\\36&-7}-1\bigg)[T^{1/5}]+\bigg(\bm{5&1\\-36&-7}-1\bigg)[T^{-1/5}]\\
    +&\bigg(\bm{5&-2\\-72&29}-1\bigg)[T^{2/5}]+\bigg(\bm{5&2\\72&29}-1\bigg)[T^{-2/5}]\bigg)\eq 0.
\end{align*}
Transferring this to $h$ gives that
\begin{align*}
    h|_k\bigg(&\bigg(\bm{5&-6\\6&-7}-1\bigg)[T^{6/5}]+\bigg(\bm{5&6\\-6&-7}-1\bigg)[T^{-6/5}]\\
    +&\bigg(\bm{5&-12\\-12&29}-1\bigg)[T^{12/5}]+\bigg(\bm{5&12\\12&29}-1\bigg)[T^{-12/5}]\bigg)\eq 0.
\end{align*}
Write
\begin{align*}
    \gamma_1\eq\m{5&-6\\6&-7},\quad\gamma_{-1}\eq\m{5&6\\-6&-7},\quad\gamma_2\eq\m{5&-12\\-12&29},\quad\gamma_{-2}\eq \m{5&12\\12&29}.
\end{align*}
It is easy to check that
\begin{align*}
 \gamma_{\pm1}&\eq T^{\pm1}ST^{\pm6}ST^{\mp1},\\
 \gamma_{\pm2}&\eq (ST^{\pm2}ST^{\mp2})^2.
\end{align*}
Since $T^6$ acts trivially on $V$ by \zcref{shape-of-v}, we have $h|_k\gamma_{\pm1}=h$. For $\gamma_{\pm2}$, using (\ref{st-matrices}), it is easy to check that
\begin{align*}
    h|_k\gamma_{\pm2}\eq h\pm(\epsilon-1)(h_1+h_2).
\end{align*}
The previous equation then simplifies to the following
\begin{align*}
    (\epsilon-1)(h_1+h_2)|_k([T^{12/5}]-[T^{-12/5}])\eq 0.
\end{align*}
Since $\epsilon=-1$, we thus obtain that
\begin{align*}
    (h_1+h_2)|_kT^{24/5}\eq h_1+h_2.
\end{align*}

Write $h_3=h_1+h_2$. Then, $h_3|_kT^{24/5}=h_3$. Also, as $\epsilon=-1$, it is easy to check that $h_3|_kS=h_3|_kT=h_3$ using (\ref{table}). Hence, we have $h_3|_kT^{-1/5}S=h_3$, where $T^{-1/5}S=\sm{-1/5&-1\\1&0}$ is elliptic of infinite order. This then implies that $h_3=0$ by \zcref{weil}. That is, $h_2=-h_1$.

Now, $h|_kT^3=-h$ and $h_1|_kT^3=-h_1$ by (\ref{table}). As $h_2=-h_1$, it then follows that $T^3$ acts as $-\id_V$ on $V$.
\end{proof}

To make use of the fact that $T^3$ acts as $-\id_V$ on $V$, we need a group-theoretic lemma.

\begin{lemma}
\zlabel{gamma6-triv}
    Let $W$ be a complex vector space with a $\PSL_2(\Z)$-action. Suppose that $T^3$ acts as $-\id_W$ on $W$. Then, $\Gamma(6)\subseteq\PSL_2(\Z)$ acts trivially on $W$.
\end{lemma}
\begin{proof}

    We first show that the normal closure $\langle\!\langle T^3\rangle\!\rangle$ of $T^3$ inside $\PSL_2(\Z)$ equals $\Gamma(3)$. Indeed,
    \begin{align*}
        \PSL_2(\Z)/\langle\!\langle T^3\rangle\!\rangle\eq\langle s,t\mid s^2=t^3=(st)^3=1\rangle\;\cong\;\PSL_2(\F_3)\;\cong\;\PSL_2(\Z)/\Gamma(3).
    \end{align*}
    Since $T^3\in\Gamma(3)$, this then implies that $\langle\!\langle T^3\rangle\!\rangle=\Gamma(3)$.

    Let $\rho:\PSL_2(\Z)\rightarrow\GL(W)$ represent the action and let $\chi$ be a character defined by
    \begin{align*}
        \chi:\PSL_2(\Z)\rightarrow\PSL_2(\F_2)\cong S_3\xrightarrow{\sgn}\{\pm1\}.
    \end{align*}
    We have $\chi(S)=\chi(T)=-1$. Let $\wt{\rho}=\chi\otimes\rho$. Then,
    \begin{align*}
        \wt{\rho}(T^3)\eq\chi(T)^3\rho(T^3)\eq(-1)(-\id_W)\eq\id_W.
    \end{align*}
    Hence, $\Gamma(3)=\langle\!\langle T^3\rangle\!\rangle\subseteq\Ker\wt{\rho}$. That is,
    \begin{align*}
        \rho(\gamma)\eq\chi(\gamma)\id_W
    \end{align*}
    for all $\gamma\in\Gamma(3)$. Since $\Gamma(2)\subseteq\Ker\chi$, this implies that
    \begin{align*}
        \Gamma(2)\cap\Gamma(3)\subseteq\Ker\rho.
    \end{align*}
    The result then follows by noting that $\Gamma(6)\subseteq\Gamma(2)\cap\Gamma(3)$.
\end{proof}

\begin{corollary}
\zlabel{V-gamma6-inv}
    The group $\Gamma(6)$ acts trivially on $V$.
\end{corollary}
\begin{proof}
    This follows directly from \zcref{t3-triv,gamma6-triv}.
\end{proof}

\begin{lemma}
\zlabel{h-gamma6-inv}
    We have $h\in S_k(\Gamma(6))$.
\end{lemma}
\begin{proof}
    It follows from \zcref{V-gamma6-inv} that $h|_k\gamma=h$ for all $\gamma\in\Gamma(6)$. Now, it suffices to note that the definition of $h$ and $f$ implies that $h$ has a Fourier expansion with only positive terms and coefficients of polynomial growth, so $h\in M_k(\Gamma(6))$ by \cite[Proposition~1.2.4]{diamond-shurman}.

    For the cuspidality, note that for any $\gamma\in\SL_2(\Z)$, we have $h|_k\gamma\in V$, so $(h|_k\gamma)(z+3)=-(h|_k\gamma)(z)$ by \zcref{t3-triv}. It then follows from \zcref{cusp-lemma} that $h|_k\gamma$ vanishes at infinity. In conclusion, $h$ vanishes at all the cusps.
\end{proof}

\begin{corollary}
    \zcref{mainthm} holds when $N=36$. That is, $f\in S_k(\Gamma_0(36))$.
\end{corollary}
\begin{proof}
    Fix $\gamma=\sm{a&b\\36c&d}\in\Gamma_0(36)$. We have
\begin{align*}
    \m{1&0\\0&1/6}\gamma\m{1&0\\0&6}\eq\m{a&6b\\6c&d}.
\end{align*}
Since $ad=36bc+1\equiv 1\pmod{6}$, we have $a\equiv d\equiv\pm1\pmod{6}$. Hence, by replacing $\gamma$ with $-\gamma$ if necessary, we always obtain that $\sm{1&0\\0&6}^{-1}\gamma\sm{1&0\\0&6}\in\Gamma(6)$, so $h|_k\sm{1&0\\0&6}^{-1}\gamma\sm{1&0\\0&6}=h$ by \zcref{h-gamma6-inv}. This then implies that $f|_k\gamma=f$ for all $\gamma\in\Gamma_0(36)$. As $f=h|_k\sm{1&0\\0&1/6}$ and $h\in S_k(\Gamma(6))$ by \zcref{h-gamma6-inv}, it follows that $f\in S_k(\Gamma_0(36))$.
\end{proof}

\bibliographystyle{amsalpha}
\bibliography{ref}

\end{document}